\documentclass[amsthm]{elsart}

\usepackage{yjsco}
\usepackage{natbib}
\usepackage{amsmath}
\usepackage{amssymb}

\numberwithin{equation}{section}

\newcommand{\qbin}[2]{\genfrac{[}{]}{0pt}{}{#1}{#2}}
\newcommand{\gp}[3]{\qbin{#1}{#2}_{#3}}
\newcommand{\Tzero}[3]{\mathrm T_0 {\textstyle(#1,#2;#3)}}
\newcommand{\Tone}[3]{\mathrm T_1 {\textstyle(#1,#2;#3)}}
\newcommand{\V}[3]{\mathrm V {\textstyle(#1,#2;#3)}}
\newcommand{\ft}{\Psi}

\begin{document}
\begin{frontmatter}

\title{Rogers-Ramanujan Computer Searches}


\author{James McLaughlin}
\address{Department of Mathematics, Room 179, 25 University Avenue \\West Chester University, West Chester PA 19383, USA}
\ead{jmclaughl@wcupa.edu}
\ead[url]{http://math.wcupa.edu/\~{}mclaughlin/}

\author{Andrew V. Sills}
\address{Department of Mathematical Sciences, 203 Georgia Avenue Room 3008,\\ Georgia Southern University, Statesboro, GA 30460-8093, USA}
\ead{ASills@GeorgiaSouthern.edu}
\ead[url]{http://math.georgiasouthern.edu/\~{}asills}

\author{Peter Zimmer}
\address{Department of Mathematics, Room 179, 25 University Avenue, \\West Chester University, West Chester PA 19383, USA}
\ead{pzimmer@wcupa.edu}

\begin{abstract}
We describe three computer searches (in PARI/GP, Maple, and Mathematica, respectively)
which led to the discovery of a number of identities of Rogers-Ramanujan type and
identities of false theta functions.
\end{abstract}


\date{2 March 2009}

\end{frontmatter}

\section{Introduction} \label{intro}

L. J.~\cite{R94} published what later came to be
known as the Rogers-Ramanujan identities:

\begin{thm}
\begin{equation} \label{RR1}
   \sum_{n=0}^\infty \frac{q^{n^2}}{(q;q)_n} =
   \frac{(q^2,q^3,q^5;q^5)_\infty}{(q;q)_\infty},
\end{equation}
and
\begin{equation} \label{RR2}
   \sum_{n=0}^\infty \frac{q^{n(n+1)}}{(q;q)_n} =
   \frac{(q,q^4,q^5;q^5)_\infty}{(q;q)_\infty},
\end{equation}
where the finite rising $q$-factorial is given by
$$(a;q)_n = \prod_{j=0}^{n-1} (1-aq^j),$$
the infinite rising $q$-factorial is given by
$$(a;q)_\infty = \prod_{j=0}^\infty (1-aq^j),$$ and
a collection of several rising $q$-factorials on the same base $q$ may
be abbreviated as
$$(a_1, a_2, \cdots a_r; q)_\infty = (a_1;q)_\infty (a_2;q)_\infty \cdots
(a_r;q)_\infty.$$
\end{thm}

Following Ramanujan, let us define
\[ f(a,b) := \sum_{n=-\infty}^\infty a^{n(n+1)/2} b^{n(n-1)/2}. \]
By Jacobi's triple product identity~\cite[p. 21, Theorem 2.8]{A76},
\begin{equation} \label{JTP}
 f(a,b) = (-a, -b, ab; ab)_\infty.
 \end{equation}

Ramanujan employed the following abbreviations for three instances
of $f(a,b)$ that occur particularly often:
\[ f(-q):= f(-q,-q^2), \qquad \varphi(q):= f(q,q), \qquad \psi(q):= f(q,q^3). \]
Notice that by~\eqref{JTP}, we have
\[ f(-q) = (q;q)_\infty, \qquad \varphi(-q):=\frac{(q;q)_\infty}{(-q;q)_\infty}, \qquad
\psi(-q) = \frac{ (q^2;q^2)_\infty }{(-q;q^2)_\infty}. \]

Rogers's paper included other identities which were similar in form in that
the left hand side is an infinite series which included rising $q$-factorials
and $q$ raised to a power that is quadratic in the summation variable, and
the right hand side is an instance of $f(a,b)$ divided by either $f(-q)$, $\phi(-q)$,
or $\psi(-q)$.
Later,~\cite{R17} published another paper
where he simplified his proof of the Rogers-Ramanujan identities and
provided additional identities of similar type.  Rogers's work was largely
ignored by the mathematical community until Ramanujan independently
rediscovered the Rogers-Ramanujan identities, and later found Rogers's
1894 paper while reading old issues of the \textit{Proceedings of the London
Mathematical Society}.

W.N.~\cite{B47,B48} undertook a careful study of Rogers's papers and
discovered the underlying engine which brought Rogers's identities
into being, and introduced the term ``Rogers-Ramanujan type
identity."  This engine was named ``Bailey's Lemma" by~\citet[p.
270]{A84}.  By calculating what~\citet[p. 26]{A86} dubbed ``Bailey
pairs," one could produce identities of the Rogers-Ramanujan type.
Bailey produced a number of such pairs in the two aforementioned
papers. As it turns out, a young Freeman Dyson was asked to referee
Bailey's papers, and in the process worked out a number of new
Rogers-Ramanujan type identities. He did not proceed in any
systematic fashion, but rather produced ``a rather haphazard
selection" which~\citet[p. 435]{B47} gladly added to his paper as he
felt ``they add[ed] considerably to the interest of the paper."

  In the late 1940's, Lucy J. Slater, then a Ph.D. student of Bailey at the
University of London, worked out many different
Bailey pairs and deduced from them a list of 130 Rogers-Ramanujan type
identities.  The now famous ``Slater list"~\citep{S52}
includes most of the Rogers-Ramanujan type identities of Rogers, Ramanujan, Bailey,
and Dyson, as well as many that were new at the time.  Although the
production of a list of 130 identities is an impressive tour-de-force,
nonetheless one is led to wonder whether some related identities
remained lurking in the background undiscovered.  Indeed this had proved
to be the case; over the years additional Rogers-Ramanujan-Slater type
identities have turned up in the work of other mathematicians; e.g.~\cite{GS83},~\cite{ S01}.

 Furthermore, since Rogers-Ramanujan type identities were first shown
to have applications to physics and Lie algebras, there has been a
virtual explosion of research in these areas.  Much of the work of
the physicists has centered around infinite families of multisum
polynomial identities as they apply to various models in statistical
mechanics.

 Closely related to identities of Rogers-Ramanujan type are identities of so-called
``false theta" functions. Noting that Ramanujan's theta series
\[ f(a,b) := \sum_{n=-\infty}^\infty a^{n(n+1)/2} b^{n(n-1)/2} = \
  \sum_{n=0}^\infty a^{n(n+1)/2} b^{n(n-1)/2} +  \sum_{n=1}^\infty a^{n(n-1)/2} b^{n(n+1)/2},  \]
let us define the corresponding \emph{false theta series} as
 \begin{equation} \label{ftdef}
 \ft(a,b) :=
  \sum_{n=0}^\infty a^{n(n+1)/2} b^{n(n-1)/2} - \sum_{n=1}^\infty a^{n(n-1)/2} b^{n(n+1)/2}.
  \end{equation}
\cite{R17} studied $q$-series expansions for many instances of
$\ft(\pm q^\alpha, \pm q^\beta)$, which he
called \emph{
false theta 
series of order $(\alpha+\beta)/2$}.

In the present paper we present the results of a number of computer searches that were undertaken with the aim of finding new Rogers-Ramanujan-type identities and new identities for false theta functions. Before coming to these searches, we briefly recall some  related investigations by other authors, in order to put the present searches in context. These other investigations also serve to show that the discovery of new series-product identities is not just interesting in its own right, but also has implications in related areas.

Firstly,~\cite{A79} called a pair of sequences $\{a_n\}$, $\{b_n\}$ a \emph{Ramanujan Pair}, if the identity
\begin{equation*}
\frac{1}{\prod_{n=1}^{\infty}(1-q^{a_n})}=1+\sum_{n=1}^{\infty}\frac{q^{b_1+b_2+\dots+b_n}}{(q;q)_n}
\end{equation*}
holds. Andrews listed four such pairs, and conjectured these were the only ones. However, two further pairs were given by~\cite{H80}, and an additional four pairs, discovered via a computer-assisted search, are given in the paper \cite{BBG86} by Blecksmith, Brillhart and Gerst. As it turns out, all of these Ramanujan pairs corresponded to existing identities of Rogers-Ramanujan  type. It would have been interesting if some of the new identities in the present paper corresponded to new Ramanujan pairs, but, alas, this is not the case.

Secondly, let the function $L(x)$ be defined by
\begin{equation}
L(x)=\sum_{n=1}^{\infty}\frac{x^n}{n^2}+\frac{1}{2} \log x \log (1-x),
\end{equation}
where the series on the right is the  \emph{dilogarithm function}.  \cite{RS81} applied asymptotic analysis to  two Rogers-Ramanujan-type identities to derive identities for the function $L(x)$. \cite{L84} extended these results to show that many other identities for $L(x)$ could be derived from known identities of Rogers-Ramanujan type. In fact, the parallels led Loxton to conjecture and prove two new series-product identities, while proving two dilogarithm identities conjectured by \cite{L82}, and also to discover a new dilogarithm identity deriving from an identity of Rogers-Ramanujan type. What is interesting in connection with finding new series-product identities experimentally is that one dilogarithm identity, conjectured by Lewin, lacks a partner on Slater's list of identities \cite{S52}, as Loxton attested, after working through all 130 of Slater's identities. Unfortunately, the present searches did not succeed in uncovering that partner either.

Thirdly, the types of computer searches described in the present paper are not the only types which may lead to the discovery of new identities of Rogers-Ramanujan-type. Indeed, we recall a sequence of papers \cite{AKK00, AK01, AKPZ01, AKPP01} in which the authors describe how to produce the $q$-Engel expansion of an infinite  $q$-series or $q$-product. Let $A= \sum\sb {n \geq \nu}c\sb n q\sp n$ be the Laurent series in $q$ corresponding to the $q$-series or $q$-product, and define $[A] = \sum\sb {\nu \leq n \leq 0} c\sb nq\sp n$. Set $a\sb 0 = [A]$, $A\sb 1 = A - a\sb 0$, and recursively define $A\sb {n+1} = a\sb n\,A\sb n - 1$, where $a\sb n = [1/A\sb n]$, $n\geq 1$. The sequence $(a\sb n)$ is the Engel sequence. The Engel expansion of $A$ is given by
\[
A = a\sb 0 + \sum\sb {n=1}\sp {\infty} \frac{1}{ a\sb 1 \cdots a\sb n}.\]
It turns out that many series - series or series-product  $q$-identities follow from the fact that a series on one side of the identity is simply the Engel expansion of the series or product on the other side of the identity. The authors use this fact to give new proofs of many existing identities. In \cite{AKPZ01}, the authors describe a \emph{Mathematica} package \textbf{Engel}, which automates the process of producing the Engel expansion of an infinite $q$-series or $q$-product, and use this package to discover a new infinite Rogers-Ramanujan type  family of
identities.

Fourthly and lastly, many partition identities have analytic representations as identities of Rogers-Ramanujan type. Thus a discovery of a new identity of one type may lead to new identity of the other type.
in \cite{A89} for example, Andrews used IBM's SCRATCHPAD to give a new proof of a partition theorem of G\"ollnitz (Andrews's result is actually a refinement of the result of G\"ollnitz), by proving a certain associated analytic identity. Other applications of computer algebra to similar investigations are described in \cite{A71} and \cite{A78}.

\vspace{10pt}

  Two of the experiments described in the present paper use PARI/GP and Maple respectively
to systematically search through Rogers-Ramanujan type series to see which
of these series have a Rogers-Ramanujan type product associated with it.
The Maple search also leant itself to a search for identities of false theta functions.
The Mathematica experiment searched for Bailey pairs of a certain form.
Of course, most of the output for all for all of the experiments revealed known identities, but
a number of the identities found appear to be new.
We list them in section~\ref{ids}. We remark that, while the false theta function identities were found through the Maple search, most  of the remaining identities turned up as output from more than one of the three searches.

   In sections \ref{PariGPSearch}--\ref{MathematicaSearch}, we describe the experimental designs.   In section~\ref{Proofs}, we indicate proofs, using a variety of methods,
of the identities presented.

\section{Some new Rogers-Ramanujan type and false theta function identities} \label{ids}
A number of the identities that have turned up in our searches fit nicely with various themes
and accordingly we have published them elsewhere~\citep{BMS08, MS08a, MS08b}.
The following additional identities were found in one or more of our searches.
In each case, we indicate a series on the left hand side equal to a ``raw product" (to show how it is
a ratio of two theta functions) as the middle member, and then provide a simplified form
of the product on the right hand side.

\begin{equation}\label{mod2-1}
  \sum_{n=0}^\infty \frac{ q^{n(n-1)} (-q;q^2)_n }{(q;q)_{2n} } = \frac{ f(1,q^2)}{ \psi(-q)} =
  (-q;q)_\infty (-1;q^2)_\infty
 \end{equation}
 \begin{equation} \label{mod2-2}
  \sum_{n=0}^\infty \frac{ q^{n(n+1)} (-q;q^2)_n }{(q;q)_{2n+1} } = \frac{ f(1,q^2)}{2\psi(-q)}
  = (-q;q)_\infty (-q^2;q^2)_\infty
 \end{equation}
 \begin{equation} \label{mod4}
 \sum_{n=0}^\infty \frac{ q^{n(n+1)} (q^2;q^2)_{n+1} }{(-q^3; q^3)_{n+1} (q;q)_n }
 =\frac{ f(-q^2,-q^2)}{ \varphi(-q^2)} = 1
 \end{equation}

 \begin{multline} \label{mod8}
  \sum_{n=0}^\infty \frac{ q^{2n(n-1)} (-q^4;q^4)_n (q;q^2)_{2n} } {(q^4;q^4)_{2n}}
  = \frac{ f(-q,-q^7)+ f(-q^3,-q^5)}{\varphi(-q^4)}\\
  = \frac{ (q,q^7; q^8)_\infty + (q^3, q^5; q^8)_\infty } { (q^4;q^8)_\infty^2}
\end{multline}

   The next identity is a partner to~\citep[Eq. (21)]{S52} and~\citep[Eq. (2.17)]{BMS08}.
 \begin{equation} \label{mod5}
    \sum_{n=0}^\infty \frac{ (-1)^n q^{n(n+2)}  (q;q^2)_n}{ (-q;q^2)_{n+1} (q^4;q^4)_n}
     = \frac{f(1,q^5)}{2\psi(-q)}
     = \frac{ (q^{10};q^{10})_\infty (q^{20};q^{20})_\infty}{(q;q^2)_\infty (q^5;q^{20})_\infty
     (q^4;q^4)_\infty}
\end{equation}

  The next identity provides an alternate series expansion for the infinite product
in~\citet[p. 154, Eq. (26)]{S52}, and is thus a partner of Slater's (22) and (26).
\begin{equation} \label{mod6-ss}
  \sum_{n=0}^\infty \frac{ q^{n^2} (-1;q)_n }{ (q;q^2)_n (q;q)_n }
  = \frac{ f(-q^3,-q^3)}{\varphi(-q)}
  =\frac{ (q^3;q^3)_\infty (q^3;q^6)_\infty }{(q;q)_\infty (q;q^2)_\infty}
\end{equation}

The next identity is a partner to~\citep[p. 154, Eq. (25)]{S52}.
\begin{equation} \label{mod6-atns}
  \sum_{n=0}^\infty \frac{q^{n(n+2)} (-q;q^2)_n}{(q^4;q^4)_n}
  = \frac{ f(-q,-q^5)}{\psi(-q)}
  = \frac{ (q^6;q^6)_\infty }{ (q^4;q^4)_\infty (q^3, q^9; q^{12})_\infty }
\end{equation}


The next identity is a partner to~\citep[p. 157, Eqs. (56) and (58)]{S52}.
\begin{equation}
\frac{1+q^3}{(1-q)(1-q^2)} + \sum_{n=1}^\infty \frac{ q^{n(n+2)} (-q;q)_{n-1} (-q;q)_{n+2}}{(q;q)_{2n+2} }
=\frac{ f(q^3,q^9)}{f(-q)}
= \frac{ (-q^3, -q^9, q^{12}; q^{12} )_\infty} { (q;q)_\infty} \label{mod12}
\end{equation}

The next identity provides an alternate series expansion for the infinite product
in~\citet[p. 159, Eq. (69)]{S52}.
\begin{multline} \label{mod16}
1+\frac{q}{(1-q)(1-q^2)}+\sum_{n=2}^\infty \frac{ q^{n^2-2} (-q^2;q^2)_{n-2} (1+q^{2n+2}) }{(q;q)_{2n} }
= \frac{f(q^2, q^{14})}{\psi(-q)}\\
= \frac{ (-q^2, -q^{14}, q^{16}; q^{16} )_\infty (-q;q^2)_\infty}  { (q^2;q^2)_\infty}
\end{multline}

The following are identities of false theta functions.
\begin{equation}
\ft(-q^3,-q) =
  \sum_{n=0}^\infty \frac{  (-1)^n q^{n(n+1)} (-q;q^2)_{n} }{ (q;q^2)_{n+1} (-q^2;q^2)_n  }
  \label{H19ft}
 \end{equation}
\begin{equation}
 \ft( q^{15}, q^{3} ) =
   \sum_{n=0}^\infty \frac{ (-1)^n q^{n(n+3)/2} (q^3;q^3)_n (1-q^{n+1})}{(q;q)_{2n+2} }
   \label{J3ft}
\end{equation}

\begin{equation}
\ft(-q^8,-q^{24}) =
  \sum_{n=0}^\infty \frac{ (-1)^{n} q^{n(n+3)/2} (q;q)_{n+1} (-q^2;q^2)_n }{ (q;q)_{2n+2}  }
  \label{K6ft}
\end{equation}
\begin{equation}
\ft(q^{22}, q^{10}) + q \ft(q^{26},q^6) =
  \sum_{n=0}^\infty \frac{ (-1)^{n} q^{n(n+3)/2} (q;q)_{n} (-q;q^2)_n }{ (q;q)_{2n+1}  }
  \label{K4ft}
\end{equation}


\section{The PARI/GP search} \label{PariGPSearch}
The pari/gp searches involved computing to high precision series of
the form
\begin{align}\label{Seq}
S&:=\sum_{n=0}^{\infty}\frac{q^{(a n^2+b
n)/2}(-1)^{cn}(d,e;q)_n}{(f,g,q;q)_n},\\
S'&:=\sum_{n=0}^{\infty}\frac{q^{(a n^2+b
n)/2}(-1)^{cn}(d,e;q^2)_n}{(f,g,q^2;q^2)_n}, \notag \\
S''&:=\sum_{n=0}^{\infty}\frac{q^{(a n^2+b
n)/2}(-1)^{cn}(d;q)_n}{(e;q^2)_{n+1}(q;q)_{n+1}}, \notag
\end{align}
for a fixed numerical value of $q$ (say $q=0.0001$), for
\[
d,e,f,g \in \{0,-1,q,-q,-q^2,q^2\}, \hspace{20pt} c \in \{0,1\},
\]
and for integers $a$ and $b$ satisfying $0 \leq |b| \leq a \leq 10$.
These choices for the forms of the series were motivated by the fact
that many series on the Slater list have similar forms.

For each particular choice of the parameters $a$, $b$, $c$, $d$,
$e$, $f$ and $g$, a numerical comparison was performed to see if
\begin{equation}\label{SPeq}
S-\prod_{j=1}^{L}(q^j;q^{L})_{\infty}^{s_j}=0,
\end{equation}
for integers $s_j$ and $L \in \{20,24,28,32,36\}$. With sufficient
precision, a small numerical value for the left side of \eqref{SPeq}
indicated either  a known identity or a potential new identity,
which then needed to be proved.

In practice, we began by creating five lists $\mathcal{L}_L$, for $L
\in \{20,24,28,32,36\}$. Here \[ \mathcal{L}_L=\{\log (q;q^L), \log
(q^2;q^L), \dots , \log (q^L;q^L)\}.
\]
We let the variables $a$, $b$, $c$, $d$, $e$, $f$ and $g$ loop
through their allowed values, and for each set of choices, computed
the series, say $S$, to high precision. For each $L$ we add $\log S$
to $\mathcal{L}_L$, to create a new list
\[ \mathcal{L}_L'=\{\log (q;q^L),
\log (q^2;q^L), \dots , \log (q^L;q^L), \log(S)\}.
\]

Next, we apply pari/gp's version of the LLL algorithm to the list
$\mathcal{L}_L'$, via the ``\emph{lindep}" command. This command
causes pari/gp to look for a zero-sum integral linear relation
amongst the elements of $\mathcal{L}_L'$, i.e. a collection of
integers $\{ b_0,b_1,\dots , b_L\}$ such that
\begin{equation}\label{logeq2}
b_0\log S + \sum_{j=1}^{L} b_j \log(q^j;q^L)_{\infty}=0.
\end{equation}
If \eqref{SPeq} holds, then such a relation will exist, in the form
\begin{equation}\label{logeq}
\log S - \sum_{j=1}^{L} s_j \log(q^j;q^L)_{\infty}=0.
\end{equation}

Pari/gp will output a set of $b_j$'s making the left side of
\eqref{logeq2} zero to within the working precision of the
\emph{lindep} command, even when no exact linear relation exist.
This usually means large absolute values for the $b_j$'s, and we
suppress this unwanted output by restricting output to cases where
$\max_{0\leq j \leq L} |b_j|<10$. This restriction meant that all
the output corresponded to either known identities, or new
identities.

Remark: It is quite likely that varying the form of the series at
\eqref{Seq}, and/or extending the ranges of the parameters $a$, $b$,
 $d$, $e$, $f$ and $g$, and/or comparing the series with combinations of $q$-products
 to  moduli $L$ other than 20, 24, 28, 32 and 36 (see \eqref{SPeq}),  may uncover yet further identities.

\section{The Maple search} \label{MapleSearch}
Our Maple search differs from the PARI/GP search in that it is uses purely
symbolic rather than numerical methods, combined with some key observations.
The reason for specifically using Maple is that we wanted to take advantage
of some procedures included in F. Garvan's Maple package
\texttt{qseries}~\citep{G99}.

First, we note that Rogers-Ramanujan type identities are generally of the
form
\[ \mbox{series} = \frac{f(a,b)  }{ \theta(q^k) } ,\]
where $a$ and $b$ are $\pm q^j$, and $j$ and $k$ are positive integers,
and where
$\theta(q)$ is one of $f(-q)$, $\varphi(-q)$, or $\psi(-q)$.
Next, note that the series expansion of $f(a,b)$ is
\[ f(a,b) = 1 + a + b + a^3b + ab^3 + a^6b^3 + a^3b^6 + a^{10} b^6 + a^6 b^{10} + \cdots \]
and thus the series expansion of $f(a,b)$ will in general be sparse, i.e. the
coefficients of most powers of $q$ will be $0$.

For example,
\[ f(-q^2,-q^3) = 1 -q^2 - q^3 + q^9 + q^{11} -q^{21} - q^{24} + \cdots. \]

Thus, our plan was to consider the series expansion
   \[  \theta(q^k) \sum_{j=0}^\infty \frac{ q^{bj^2 + cj} \prod_{i=1}^h (n_{i_1} q^{n_{i_2}} ; q^{n_{i_3}})_{j+n_{i_4}} }
   { (q^m;q^m)_j  \prod_{s=1}^r ( d_{s_1} q^{s_2} ; q^{s_3}  )_{j+d_{s_4} }}
   \]
   where $1\leqq m \leqq 8$, $1\leqq b\leqq 4$, $0 \leqq c \leqq 4$, $1\leqq k \leqq 4$,
   $n_{i_1} = \pm 1$,
   $\lceil \frac{n_{i_1}}{2} \rceil \leqq {n_{i_2}} \leqq 6$, $n_{i_2} \leqq n_{i_3} \leqq 6$,
   $0 \leqq n_{i_4} \leqq 1$,
     $d_{i_1} = \pm 1$,
   $\lceil \frac{d_{i_1}}{2} \rceil \leqq {d_{i_2}} \leqq 6$, $d_{i_2} \leqq d_{i_3} \leqq 6$,
   $0 \leqq d_{i_4} \leqq 1$,  $0\leqq h \leqq 3$, $0 \leqq r \leqq 2$,
and check to see whether the Taylor series expansion was sparse.  Our sparseness criterion
was simply to check whether the Taylor expansion, truncated at the $q^{55}$ term, had
less than $12$ terms.  If so, then we used procedures from F. Garvan's \texttt{qseries}
Maple package~\citep{G99} to try to factor the Taylor series into a product, and, if possible, identify
the series as an instance of $f(a,b)$.
On the other hand, if the series expansion of
\[ \sum_{j=0}^\infty \frac{ q^{bj^2 + cj} \prod_{i=1}^h (n_{i_1} q^{n_{i_2}} ; q^{n_{i_3}})_{j+n_{i_4}} }
   { (q^m;q^m)_j  \prod_{s=1}^r ( d_{s_1} q^{s_2} ; q^{s_3}  )_{j+d_{s_4} }}
   \] was found to be sparse, but not expressible as an infinite product, it was checked to see
 if it was a false theta function.

\section{The Mathematica search} \label{MathematicaSearch}
For this search, we used Mathematica because we wanted to take
advantage of A. Riese's excellent implementation of the $q$-analog
of Zeilberger's algorithm, which is available as the Mathematica
package \texttt{qZeil.m}~\citep{PR97}.

A classical method for establishing Rogers-Ramanujan identities is via
inserting {Bailey pairs} into limiting cases of Bailey's lemma.

  A pair of sequences $\Big(  \{ \alpha_n (x,q) \}_{n=0}^\infty ,   \{ \beta_n (x,q) \}_{n=0}^\infty \Big)$
is called a~\emph{Bailey pair relative to $x$} if
  \begin{equation} \label{BPdef}
   \beta_n(x,q) = \sum_{r=0}^n \frac{ \alpha_r (x,q)}{ (q;q)_{n-r} (xq;q)_{n+r} }
   \end{equation}
for all nonnegative integers $n$.

Furthermore, Bailey's lemma~\citep{A86} implies that if   $( \alpha_n (x,q), \beta_n (x,q) )$
form a Bailey pair,
then
\begin{equation} \label{WBL}
  \sum_{n=0}^\infty x^n q^{n^2} \beta_n (x,q)
  = \frac{1}{(xq;q)_\infty} \sum_{n=0}^\infty x^n q^{n^2} \alpha_n (x,q),
  \end{equation}
\begin{equation} \label{ATNSBL}
  \sum_{n=0}^\infty x^n q^{ n^2} (-q;q^2)_n \beta_n (x,q^2)
  = \frac{(-xq; q)_\infty}{(xq^2;q^2)_\infty} \sum_{n=0}^\infty
  \frac{x^n q^{n^2} (-q;q^2)_n}{(-xq;q^2)_n} \alpha_n (x,q^2),
  \end{equation}
  \begin{equation} \label{ATNSnegBL}
  \sum_{n=0}^\infty (-1)^n x^n q^{ n^2} (q;q^2)_n \beta_n (x,q^2)
  = \frac{(xq; q)_\infty}{(xq^2;q^2)_\infty} \sum_{n=0}^\infty
  \frac{(-1)^n x^n q^{n^2} (q;q^2)_n}{(xq;q^2)_n} \alpha_n (x,q^2),
  \end{equation}
\begin{equation} \label{SSBL1}
   \sum_{n=0}^\infty q^{n(n+1)/2} (-1;q)_n \beta_n(1,q) = \frac{2}{\varphi(-q)}
   \sum_{n=0}^\infty \frac{  q^{n(n+1)/2}}{1+q^n} \alpha_n(1,q),
 \end{equation}
 and
\begin{equation}\frac{1}{1-q}\sum_{n=0}^\infty  (-1)^n q^{n(n+1)/2} (q;q)_n \beta_n(q,q)
=  \sum_{r=0}^\infty  (-1)^r q^{r(r+1)/2}  \alpha_r (q,q). \label{FBL}
  \end{equation}

See \citet[page 5]{McLSZ08} for these transformations.
The trick is to find Bailey pairs in which the $\alpha_n$, when inserted into the right
hand side of~\eqref{WBL},~\eqref{ATNSBL}, or  \eqref{SSBL1},
with $x$ set to $1$ or some power of $q$, produces an instance of $f(a,b)$, while
simultaneously the $\beta_n$ is expressible as a finite product.

  In our search, we chose $\alpha_n$'s to conform to the aforementioned requirement,
and inserted them into the right hand side of~\eqref{BPdef}, then used the summand
of the series
as input for
the \texttt{qZeil} function in Riese's \texttt{qZeil.m} Mathematica package.
If a first order recurrence was found, the result was recorded, and could then be
iterated to find an expression for $\beta_n$ as a finite product.

For example, let us choose $\alpha_n$ so that $\alpha_{0} = 1$, $\alpha_{2m+1} = 0$,
and $\alpha_{2m} = q^{2m^2 - 3m} (1+q^{6m}). $
Then the right hand side of~\eqref{WBL}, with the aid of~\eqref{JTP}, can be seen to
be equal to $f(q^3,q^9)/f(-q) = (-q^3, -q^9, q^{12}; q^{12})_\infty / (q;q)_\infty$.  On the other hand,
the $q$-Zeilberger algorithm reveals that
\begin{equation*} \beta_n = \sum_{r=0}^n \frac{\alpha_r}{(q;q)_{n+r} (q;q)_{n-r} } \end{equation*}
satisfies the recurrence
\begin{equation} \label{recur}
\beta_n = \frac{ (1+q^{n-2})(1+q^{n+1}) }{ (1-q^{2n})(1-q^{2n-1}) } \beta_{n-1}, \text{ for }n\geq 3.
\end{equation}
Equation~\eqref{recur}, together with the initial conditions calculated from~\eqref{BPdef}
imply
that
\begin{align*} 
\beta_0 &= 1,\\
\beta_1 &= \frac{1}{(1-q)^2} \\
\beta_2 &= \frac{ (1+q^{-1})(1+q^2)(1+q^3)}{(1-q)(1-q^2)(1-q^3)(1-q^4)}\\
\beta_n &= \frac{ (-q^{-1};q)_n (-q^2;q)_n }{ 2 (q;q)_{2n} } \qquad \mbox{if $n>2$}.
\end{align*}

Thus, inserting this Bailey pair into~\eqref{WBL} yields identity~\eqref{mod12}.
Inserting this same Bailey pair into~\eqref{ATNSBL} yields identity~\eqref{mod16}.

\section{Proofs of the identities} \label{Proofs}
Identities~\eqref{mod2-1} and~\eqref{mod2-2} each follow
from inserting a Bailey pair of~\citet[p. 468, F(3) and F(4) resp.]{S51} into~\eqref{ATNSBL}.
It was pointed out by one of the referees that these identities also follow as special  cases
($a=-q$, $b=q$ and $a=-q$, $b=q^3$ respectively, after replacing $q$ with $q^2$) of the following limiting form of Heine's $q$-Gauss summation:
\[
\sum_{n=0}^{\infty}\frac{(-1)^nq^{n(n-1)/2}(a;q)_n}{(b,q;q)_n}\left(\frac{b}{a}\right)^n
=\frac{(b/a;q)_{\infty}}{(b;q)_{\infty}}.
\]

Since they were found and proven by different methods, we include them. They also serve as a reminder that it is necessary to be careful when discovering ``new" identities, as many special cases of Heine's $q$-Gauss and other similar identities give identities of Rogers-Ramanujan type.  Indeed Slater's list  contains three such special cases of the $q$-Gauss identity: ~\citep[Eqs. (4), (47), (51)]{S52}.   

Identity~\eqref{mod4} is actually a limiting ($N\to\infty$) case of the polynomial identity
\begin{multline}\label{mod4poly}
 \sum_{i,j,k\geqq 0} (-1)^{i+k} q^{j^2+j+i^2+i+3k} \gp{j+1}{i}{q^2} \gp{j+k}{k}{q^3} \gp{N-2i-j-3k}{j}{q}
 \\= \sum_{j\in\mathbb{Z}} (-1)^j q^{2j^2} \Tone{N}{2j}{q},
\end{multline}
where the $q$-binomial co\"efficient is defined
\begin{equation*}
 \gp{A}{B}{q} := \left\{ \begin{array}{ll}
   \displaystyle{  \frac{(q;q)_A}{(q;q)_B (q;q)_{A-B}} }&\mbox{ if $0\leqq B \leqq A$,}\\
     0 &\mbox{otherwise,}
     \end{array} \right.
\end{equation*}
and
\begin{equation*}
 \Tone{n}{a}{q} := \sum_{j=0}^m (-1)^j q^j \gp{m}{j}{q^2} \gp{2m-2j}{m-a-j}{q}
\end{equation*}
is one of several $q$-analogs of the trinomial co\"efficient first defined by~\citet[p. 299, Eq. (2.9)]{AB87}.

  Starting with~\eqref{mod4}, both sides of~\eqref{mod4poly} can be conjectured and
proved using the methods described in detail by the second author~\citep{S04}.

 Similarly, Identity~\eqref{mod8} is a limiting ($N\to\infty$) case of the polynomial identity
 \begin{multline}\label{mod8poly}
 1+\sum_{i,j,k\geqq 0} (-1)^i q^{2j^2-2j+i^2+4k} \gp{2j}{i}{q^2} \gp{j+k-1}{k}{q^8} \gp{N-i-2k}{j}{q^4}\\
 = \sum_{j\in\mathbb Z} (-1)^j q^{4j^2+j}(1+q^{2j})\V{N}{2j+1}{q^2},
 \end{multline}
 where the function
   \[ \V{m}{a}{q} := \Tone{m-1}{a}{q} + q^{m-a}\Tzero{m-1}{a-1}{q} \]
was introduced by the second author in~\cite[p. 7, Eq. (1.23)]{S03} and
   \[ \Tzero{m}{a}{q} := \sum_{j=0}^m (-1)^j \gp{m}{j}{q^2} \gp{2m-2j}{m-a-j}{q} \]
is another $q$-trinomial co\"efficient found by~\citet[p. 299, Eq. (2.8)]{AB87}.

Identity~\eqref{mod5} follows from inserting a Bailey pair of~\citet[p. 469, G(2)]{S51}
into~\eqref{ATNSnegBL}, together with the observation that
$ 2 f(q^5, q^{10}) = f(1,q^5)$ by~\citep[p. 220, Eq. (3.2)]{B51}.
Slater's Bailey pair G(2) appears in~\citep{S51} with a misprint; the correct
form is given by
\[ \beta_n(q,q) = \frac{1}{(q^2;q^2)_n (-q^{3/2};q)_{n} } \]
and {\allowdisplaybreaks
\[ \alpha_n(q,q) = \left\{ \begin{array}{ll}
   1 &\mbox{if $n=0$}\\
   -q^{3r^2 + \frac 72 r + 1} \left( \frac{1-q^{2r+\frac 32}}{1-q^{1/2}} \right) & \mbox{if $n=2r+1$} \\
    q^{3r^2 + \frac 12 r }\left( \frac{1-q^{2r+\frac 12}}{1-q^{1/2}} \right) & \mbox{if $n=2r>0$}.
    \end{array} \right.
    \]
    }

Identity~\eqref{mod6-ss} follows from inserting a Bailey pair of~\citet[p. 469, C(5)]{S51}
into~\eqref{SSBL1}.

Identity~\eqref{mod6-atns} follows from inserting a Bailey pair of~\citet[p. 469, E(4)]{S51}
into~\eqref{ATNSBL}.

Identities~\eqref{mod12} and~\eqref{mod16} follow from
a Bailey pair that is not in Slater's paper.  As discussed in
the previous section, the Bailey pair was found using the $q$-Zeilberger algorithm,
therefore the proof (via rational function proof certificate) was produced by the
\texttt{qZeil.m} Mathematica package.  As $q$-WZ certification has been explored extensively
in the literature, we refer the interested reader to~\cite{ET90, K93, P94, PR97, WZ90, WZ92, Z90, Z91}
and simply reveal that by the rational function certifying the recurrence~\eqref{recur} is
\begin{multline*}-\frac{q^{n-2 s-1} \left(q^{2 s}-q^n\right) \left(q^{2
   s+1}-1\right) \left(q^{2 s+1}-q^n\right)
  }{\left(q^n-1\right)
   \left(q^n+1\right) \left(q^n-q\right)
   \left(q^n-q^2\right) \left(q^{2 n}-q\right) \left(q^{2
   s}+1\right) \left(-q^{2 s}+q^{4 s}+1\right)} \\
   \times  \left(q^n-q^{n+2}+q^{2 n+1}+q^{2 s+1}-q^{2 s+2}
    +q^{2s+4}+q^{n+2 s+1}-q^{n+2 s+3}-q^{2 n+2 s} \right. \\
    \left. +q^{2 n+2s+2}-q^{2 n+2 s+3}-q^{4 s+3}+q^{n+4 s+2}-q^{n+4
   s+4}+q^{2 n+4 s+3}-q\right)
 \end{multline*}
and appears as the result of issuing the~\texttt{Cert[]} function call
after finding the appropriate recurrence via the following call of
the \texttt{qZeil} function.
\begin{verbatim}
qZeil[ q^(2s^2 - 3s)(1 + q^(6s))/2
      / qPochhammer[q, q, n - 2s]/qPochhammer[q, q, n + 2s],
      {s, -Infinity, Infinity}, n]
\end{verbatim}
  Notice that Paule's creative symmetrization~\citep{P94} has been used.

   The identities~\eqref{H19ft} and~\eqref{J3ft} follow from inserting Bailey pairs
in~\citep[H(19) and J(3) resp.]{S52} into~\eqref{FBL}, while identities~\eqref{K6ft} and~\eqref{K4ft}
follow from inserting the Bailey pairs in~\citep[p. 471, 6th and 4th entries resp. in the table]{S51}
into~\eqref{FBL}.

\section{Conclusion}
 Our aim here has been to show that computer algebra systems can be used effectively to search for identities of Rogers-Ramanujan type. The searches described here are quite straightforward in nature.  They could certainly be
generalized and allowed to run for longer periods of CPU time.  The truly time consuming part
of the exercise is to sift through large amount of output in the hopes of finding identities that
have not previously appeared in the literature.


\bibliographystyle{elsart-harv}


\end{document}